\documentclass{amsart}
\usepackage{amsfonts,amssymb,amsmath,amsthm, graphicx, color, tikz,tikz-cd,soul}
\usepackage{url}
\usepackage{enumerate}
\usepackage[normalem]{ulem}
\usepackage{enumerate}

\newtheorem{theorem}{Theorem}[section]
\newtheorem{lemma}[theorem]{Lemma}
\newtheorem{proposition}[theorem]{Proposition}
\newtheorem{corollary}[theorem]{Corollary}
\theoremstyle{definition}
\newtheorem{definition}[theorem]{Definition}
\newtheorem{rem}[theorem]{Remark}
\newtheorem{example}[theorem]{Example}
\numberwithin{equation}{section}

\newcommand{\bte}{\begin{theorem}\quad  }
\newcommand{\ete}{\end{theorem} }
\newcommand{\bpr}{\begin{proposition}\quad  }
\newcommand{\epr}{\end{proposition} }
\newcommand{\ble}{\begin{lemma}\quad }
\newcommand{\ele}{\end{lemma}}
\newcommand{\bco}{\begin{corollary}\quad }
\newcommand{\eco}{\end{corollary} }
\newcommand{\bex}{\begin{example}\quad \rm }
\newcommand{\eex}{\end{example} }
\newcommand{\bde}{\begin{defi}\quad \rm }
\newcommand{\ede}{\end{defi} }
\newcommand{\brm}{\begin{rem} \quad \rm}
\newcommand{\erm}{\end{rem} }
\newcommand{\bpf}{\begin{proof}[\bf{Proof.\quad}] \rm}
\newcommand{\epf}{ \end{proof}}
\newcommand{\bdm}{\begin{displaymath} }
\newcommand{\edm}{\end{displaymath} }
\newcommand{\be}{\begin{eqnarray*}}
\newcommand{\ee}{\end{eqnarray*}  }

\newcommand{\lb}{\label}

\makeatletter
\newcommand\cupop{\mathop{\operator@font \cup}\nolimits}
\makeatother
\numberwithin{equation}{section}

\author{Mohammad Roueentan}
\address{Department of Mathematics, College of Engneering, Lamerd Higher Education Center, Lamerd, Iran.}
\email{rooeintan@lamerdhec.ac.ir}

\author{Mojtaba Sedaghatjoo$^{^*}$}
\address{Department of Mathematics, College of Sciences, Persian Gulf University, Bushehr, Iran.}
\email{Sedaghat@pgu.ac.ir}
\thanks{$^{^*}$Corresponding author}

\keywords{uniform, irreducible, subdirectly irreducible, large.}
\subjclass[2010]{20M50, 20M30.}
\begin{document}
\title{On Uniform Acts over Semigroups}
\begin{abstract}
The paper is devoted to the investigation of uniform notion for acts over semigroups perceived as an overclass of subdirectly irreducible acts. We establish conditions to fill the gap between these classes of acts. Besides we prove that uniform acts with two zeros are subdirectly irreducible.  Ultimately we investigate monoids which are uniform as right acts over themselves and we characterize regular ones.
\end{abstract}
\maketitle
\section{INTRODUCTION AND PRELIMINARIES}
The study of uniform notion for rings and modules was initiated by A. W. Goldie, (See \cite{goldie1,goldie2})  which led into the uniform dimension of modules and consequently, homological classification of rings with regard to uniform modules. On the other hand, due to the embedding of any algebra in a subdirect product of subdirectly irreducible algebras of the same type (Birkhoff's theorem), uniform modules as an over-class of subdirectly irreducible modules were a prominent issue for investigation. Hereby, the counterpart notions for acts over semigroups are deserved to be considered which construct the main body of this work. Apparently, the first investigation on subdirectly irreducible acts goes back to 1974, by E.N. Roiz and was followed by some authors for instance Khozhukhov in \cite{kozh1,kozh2,kozh3,kozh}. Evidently, a pioneering work on uniform acts was initiated by Feller and Gantos in \cite{feller} which they characterized injective uniform acts as the ones with a local endomorphism monoid of their injective envelops. In this paper we investigate uniform acts with a fundamental account on the basic notions and their interrelationships with injective  and subdirectly irreducible acts. We prove that uniform acts possessing two zeros are subdirectly irreducible and characterize subdirectly irreducible acts over right zero semigroups with a different approach than what is done in \cite{kozh} (in Russian). It is known that for a uniform module $M$, an endomorphism is monomorphism if and only if it is not nilpotent. We prove that the same is valid for uniform acts. Ultimately we clarify and in some cases characterize some classes of uniform semigroups as right acts over themselves.

Throughout this paper, $S$ will denote a semigroup. To every semigroup $S$ we can associate the monoid $S^1$ with the identity element $1$ adjoined if necessary. Indeed, $S^1=\begin{cases} S & \text{if} ~S~ \text{has an identity element,} \\ S  \cup\{1\} & \text{Otherwise}. \end{cases}$ We denote by $A\sqcup B$ the disjoint union of  sets $A$ and $B$. A right $S$-act $A_S$ (or an act $A$, if there is no danger of ambiguity) is a non-empty set together with an action $ \mu: A\times S \longrightarrow A ,\,as:=\mu(a,s)$, such that $a(st)=( as)t$ and in the case that $S$ is a monoid with the identity element $1$, $a.1=a$ for each $a\in A$ and $s,t \in S$. An element $\theta \in A$ is said to be a zero element if $\theta s=\theta$ for all $s\in S$. Moreover, the one element act, denoted by $\Theta=\{\theta\}$, is called the zero act. Recall that an act is called simple ($\theta$-simple) if it contains no (non-zero) subacts other than itself. Also recall that for an element $a$ of an $S$-act $A$, by $\lambda_a$ we denote the homomorphism from $S_S$ to $A$ defined by $\lambda_{a}(s)=as$ for every $s\in S$. An equivalence relation  $\rho$ on an $S$-act $A$ is called a right congruence on $A$ if  $a~\rho ~a'$ implies $(as)~ \rho ~(a's)$ for every $a,a' \in A, s \in S$. For an act $A$ the diagonal relation $\{(a,a)\,|\,a\in A\}$ on $A$ is a congruence on $A$ which is denoted by $\Delta_A$. Also if $B$ is a subact of $A$, then the congruence $(B\times B)\cup \Delta_A$ on $A$ is denoted by $\rho_B$ and is called Rees congruence by the subact $B$. Also the set of all congruences on $A$ is denoted by Con($A$). For $a,b\in A, $ the monocyclic congruence on $A$ generated by the pair $(a,b)$ is denoted by $\rho(a,b)$ which for $x,y \in A, x \rho(a,b)y$ if and only if $x=y$ or there exist $p_1,p_2,\ldots ,p_n,q_1,q_2,\ldots ,q_n\in A, w_1,w_2,\ldots ,w_n\in S^1$ where for every $i=1,2,\ldots ,n, (p_i,q_i) \in \{(a,b), (b,a)\}$, with
\begin{alignat*}{3}
x =&p_1 w_1 & &q_2 w_2 =p _3
w_3 & &~\cdots ~q_n w_n =y.\\
&q_1 w_1 \,&= \,\,&p_2  w_2 & &~\cdots ~
\end{alignat*}

We recall that an $S$-act $A$ is called subdirectly irreducible if every set of congruences $\{\rho_i \,| \,i\in I \}$ on $A$ with $\bigcap_{_{i\in I}} {\rho_i} = \Delta_A$, contains $\Delta_A$, indeed, the set of nondiagonal congruences has a least element. Also an act $A$ is called irreducible if for any two congruences $\rho$ and $\lambda$ on $A$, $\rho \bigcap\lambda =\Delta_A$ implies that $\rho=\Delta_A$ or $\lambda=\Delta_A$, equivalently, intersection of finite nondiagonal congruences is nondiagonal. It is clear that a finite $S$-act is subdirectly irreducible if and only if it is irreducible.

Considering the unique decomposition of acts into indecomposable acts over monoids, it can be routinely observed that the same structure is valid for acts over semigroups. Indeed, for a semigroup $S$, any right $S$-act $A_S$ has a unique decomposition into indecomposable acts which indecomposable components of $A_S$ are the equivalence classes of the relation $\sim $ on $A_S$  defined in \cite{Ren} by $a\sim b$ if there exist $s_1,s_2, \ldots ,s_n,t_1,t_2, \ldots ,t_n \in S^1,~a_1,a_2, \ldots , a_{n-1} \in A_S$ such that
\[ as_1=a_1t_1,~ a_1s_2=a_2t_2,~ a_2s_3=a_3t_3, \ldots,a_{n-1}s_n=bt_n. \] Based on the second approach to indecomposable components of an act, if $S$ is a group, then indecomposable acts are the same simple acts. For a thorough account on the preliminaries, the reader is referred to \cite{kilp}.

\section{UNIFORM ACTS}

In this section we bring out preliminary and basic properties of uniform acts.  First recall that a subact $B$ of an act $A$ is called large in $A$ (or $A$ is called an  essential extension of $B$), denoted by $B \subseteq' A$, if any $S$-homomorphism $g: A \longrightarrow C$ such that $g|_B$ is a monomorphism is itself a monomorphism. One may routinely observe that a subact $B$ of an act $A$ is large in $A$ if and only if for every right congruence $\rho \in$ Con$(A)$, $\rho_B \cap \rho = \Delta_A$ implies that $\rho = \Delta_A$. It is known that for acts $A\subseteq B \subseteq C$, $A \subseteq' C$ if and only if $A \subseteq' B$, $B \subseteq' C$( see \cite[Lemma 3.16]{kilp}). We know that for a ring $R$, a submodule $K$ of an $R$-module $M$ is large in $M$ if and only if every non-zero submodule of $M$ has non-zero intersection with $K$, nevertheless, for acts over semigroups this condition is only necessary. Indeed, for a semigroup $S$, if a subact $B$ of an act $A$ is large in $A$ and $X$ is a non-zero subact, then due to the canonical homomorphism from $A$ to $A/{\rho_X}$, $B$ has non-zero intersection with $X$. To observe that the mentioned condition is not sufficient, take $S=\{1,e,f\}$ where $e,f$ are right zero elements and $1$ is the identity. Then the right ideal $I=\{e,f\}$ is the only proper right ideal of $S$ which is a retract of $S$ and hence is not large in $S$. Thus $I$ satisfies the condition but it is not large in $S$. Herein, for an act $A$, the injective envelope of $A$, indeed, the least injective extension of $A$ (equivalently the maximal essential extension of $A$) is denoted by $E(A)$. Now it can be readily verified that  an act $B$ is large in $A$ if and only if $E(A)=E(B)$.

\begin{definition} For a semigroup $S$, a non-zero $S$-act $A$ is called $uniform$ if every non-zero subact is large in $A$. Also a semigroup $S$ is called right (left) uniform if the right (left) $S$-act $S_S$ ($_SS$) is uniform.
\end{definition}
Throughout this paper, unless otherwise stated, the term ``uniform semigroup" stands for ``right uniform semigroup". Readily, for a semigroup $S$, right $S$-acts with at most two elements are uniform which we shall perceive them as trivial uniform acts.
Regarding the above arguments concerning subdirectly irreducible acts and large subacts the following results are obtained.

\bco \lb{co10}
Over a monoid $S$ a right $S$-act $A$ is uniform if and only if $E(A)=E(B)$ for every non-zero subact $B$ of $A$.
\eco
\bco \lb{co2} Let $S$ be a semigroup and $A$ be a non-zero uniform $S$-act. If $B,C$ are non-zero subacts of $A$, then $|B\cap C|\geq 2$.
\eco

 For an $S$-act $A$, by $Z(A)$, we mean the zero elements of $A$. The next corollary is a straightforward established property of uniform acts.
\bco \lb{co18} Let $A$ be a uniform act over a semigroup $S$ and $B$ be a non-zero subact. Then either $Z(A)\subseteq B$ or $Z(A)\cap B=\emptyset$.
\eco
\bpf Suppose that $\theta_1\in Z(A)\cap B$ and $\theta_2 \in Z(A)\backslash B$. Thus for non-zero subacts $\{\theta_1,\theta_2\}$ and $B$ we have $|\{\theta_1,\theta_2\}\cap B|=1$ which leads to a contradiction.
\epf
\bco \lb{co5}
Let $A$ be a uniform act over a semigroup $S$. Then $|Z(A)|\leq 2$.
\eco
\bpf
Suppose by way of contradiction that $|Z(A)|> 2$. Let $a,b,c,$ be three distinct elements of $Z(A)$. If we take $B= \{a,b\}$, then trivially $B$ is not large in $Z(A)$ which is a contradiction in view of Corollary \ref{co2}.
\epf
By virtue of corollaries \ref{co2}, \ref{co18} and \ref{co5}, for any uniform act with two zeros, $Z(A) \subseteq B$ for any nonzero subact $B$ of $A$. The above corollary leads us to determine the general structure of decomposable uniform acts.
\bpr \lb{pr9}
Suppose that $S$ is a semigroup and $A$ is a decomposable uniform $S$-act. Then either $A=\Theta\amalg \Theta$ or $A= B\amalg \Theta$, where $B$ is an indecomposable uniform act with no zero element.
\epr
\bpf
 If $A=A_1\amalg A_2$ is a decomposition of $A$, then Corollary \ref{co2} necessitates that for some $i\in \{ 1,2\}$, $A_i=\Theta$ and hence $A=\Theta\amalg \Theta$ or $A= B\amalg \Theta$, where in light of Corollary \ref{co18} $B\cap Z(A)= \emptyset$. It is implicitly concluded that decomposable uniform acts contain zero elements. Since $B$ does not possesses zero element we conclude that $B$ is indecomposable.
\epf
Now we conclude that over a left zero semigroup the only decomposable uniform act is $\Theta\amalg \Theta$. On the other hand if $A$ is a non-zero act over a group, since any indecomposable act is simple, by virtue of Proposition \ref{pr9}, $A$ is uniform if an only if $A$ has one of the following structures:
\begin{enumerate} [{\rm i)}]
\item $A= \Theta\amalg \Theta$.
\item $A= B\amalg \Theta$, where $B$ is a  simple act.
\item $A$ is a non-zero simple act.
\end{enumerate}

%\bpf
%The necessity part follows by the previous theorem. For the sufficiency part suppose $B$ is a non-zero subact of $A$ and note that $B\subseteq A\subseteq E(A)=E(B)$. Now since $B\subseteq'  E(B)$, $B \subseteq' A$.
%\epf
Considering the fact that a subact $B$ of  an act $A$ is large if and only if for every right congruence $\rho \in$ Con$(A)$, $\rho_B \cap \rho = \Delta_A$ implies that $\rho = \Delta_A$, the following proposition is followed.
\bpr \lb{pr1}
For a semigroup $S$ every (subdirectly) irreducible $S$-act is uniform.
  \epr

 \bex \lb{ex1}
Let $G$ be a group. Considering $G$ as a right $G$-act, we can designate the right congruence $\rho_H=\{(x,y)\,|\, xy^{-1} \in H\}$,  for any subgroup $H$ of $G$. This assignment provides a one to one order preserving correspondence between the right congruences on $G$ and its subgroups. Note that the converse of the above proposition is not true, though, Corollary \ref{co8}, indicates that uniform acts with two  zeros are subdirectly irreducible. For instance let $S$ be an abelian group with no minimal subgroup, then clearly the right $S$-act $S_S$ is uniform but it is not subdirectly irreducible.
\eex

\brm  Note that for any  non-zero semigroup $S$, obviously $S$ is not large in $S\sqcup S$ and diagonal subact $\Delta_S=\{(x,x)\,|\,x\in S\}$  is not large in $(S\times S)_S$ for projection on one variable from $S\times S$ to $S$. Hereby, uniform notion is not preserved under products and coproducts.
\erm

It should be mention that for any semigroup $S$, generally, decomposable acts are not uniform. thereby, in what follows we characterize semigroups over which all indecomposable acts are uniform. First we recall some ingredients.

Consider monomorphisms $j_i:U\longrightarrow X_i$ for $i \in\{1,2\}$ in a category {\bf \emph{C}}. The pushout $((q_1,q_2),Q)$ of $j_1$ and $j_2$ is called an amalgam of $X_1$ and $X_2$ by the common subobject $U$ and $Q$ is denoted by $X_1\coprod ^U X_2$.

Analogously defined multiple amalgams are denoted by
$$ \coprod_{i\in I}^U X_i \text{~~for~~} U,X_i\in  {\bf \emph{C}}, i\in I. $$

  In {\bf Act-}$S$ amalgams exist and as a canonical amalgam we can have $Q=(X_1\sqcup X_2)/\nu$ where $X_1\sqcup X_2$ is the disjoint union of $X_1$ and $X_2$ and $\nu$ is the congruence relation on  $X_1\sqcup X_2$ generated by all pairs $(j_1(u),j_2(u)),u\in U$.

\bpr \lb{pr6}
Let $S$ be a semigroup. All indecomposable $S$-acts are uniform if and only if $S$ is a group. \epr

\bpf {\bf Necessity.} Let $s\in S$ and $sS\neq S$. Then the canonical epimorphism from $S \coprod ^{sS}S$ to $S$ is not a monomorphism but it is monomorphism on both copies of $S$. Thus the indecomposable act $S \coprod ^{sS}S$ is not uniform, a contradiction. Then $S$ is right simple. On the other hand if $S$ does not have an identity element, $S_S$ is a proper subact of he indecomposable act $S^1_S$. Analogously, the indecomposable acts $S^1 \coprod ^{S}S^1$ leads to a contradiction. Therefore $S$ is a group.
Regarding the fact that indecomposable acts on groups are simple, the sufficiency part follows.
\epf
 \ble \lb{le3}
 Let $S$ be a semigroup and $B \subseteq A$ be right $S$-acts. Then the following assertions are established.
\begin{enumerate} [{\rm i)}]
\item  If $A$ is  subdirectly irreducible (uniform), then $B$ is subdirectly irreducible (uniform).
\item  If $B$ is a large subdirectly irreducible (uniform) subact of $A$, then $A$ is subdirectly irreducible (uniform).
    \end{enumerate}
  \ele
\bpf
The proof for uniform acts is obvious and thus we only prove the results for subdirectly irreducible acts.
 \\
(i) Suppose that $A$ is subdirectly irreducible and $\{ \rho_i\,|\, i\in I\} \subseteq$ Con($B$) such that $\bigcap\limits _{_{i\in I}} {\rho_i} = \Delta_B$. Thus $(\bigcap\limits _{_{i\in I}} {\rho_i})\cup \Delta_A = \Delta_A$ and so $\bigcap\limits_{_{i\in I}}({\rho_i}\cup \Delta_ A) = \Delta_A$. Now by assumption ${\rho_j}\cup \Delta_A = \Delta_A$ for some $j\in I$ and the result follows.\\
 (ii) Suppose that $B$ is a large subdirectly irreducible subact in $A$. If $\{ \rho_i\,|\, i\in I\} \subseteq$ Con($A$) with $\bigcap\limits_{_{i\in I}} {\rho_i} = \Delta_A$, then  $\bigcap\limits_{_{i\in I}} ({\rho_i}\cap (B\times B)) = \Delta_B$. Now by assumption, ${\rho_j}\cap (B\times B)=\Delta_B $ for some $j\in I $ and since $B \subseteq' A$, by \cite[Lemma 3.1.15]{kilp}, ${\rho_j}= \Delta_A$.
\epf

\bpr \lb{pr4}
Let $S$ be a semigroup and $A$ be an $S$-act which has no zero element. Then $A$ is subdirectly irreducible (uniform) if and only if $ A\amalg \theta$ is subdirectly irreducible (uniform).
\epr
\bpf
By Lemma \ref{le3}, it is sufficient to show that $A$ is a large subact of $A \amalg \theta$. Suppose that $f:A \amalg \theta\longrightarrow C$ is a homomorphism and $f\mid_A$ is a monomorphism. Let $f(a)=f(b)$. If $f(a)=f(\theta)$ for some $a\in A$ then for every $s \in S$, $f(as)=f(\theta s)=f(\theta)=f(a)$ and then $as=a$ for every $s\in S$. Therefore $a$ is a zero element of $A$, a contradiction.
\epf
\brm \lb{rm3}
In light of Lemma  \ref{le3}, envelop of uniform acts are uniform acts which are indeed injective and maximal in the class of uniform acts. On the other hand any maximal uniform act in the class of uniform acts is an injective uniform act (note that injective acts are the ones with no proper essential extension). Thereby, from this aspect, any uniform act is contained in a maximal uniform act or equivalently in an injective uniform act.
%Now, for  uniform acts do not have any proper injective subact, we obtain the followin theorem.
\erm
%I will prove the following later. The last part I paused.
%\bte \lb{th4}
%Let $A$ be an act over a semigroup $S$. The following are equivalent:
%\begin{enumerate} [{\rm i)}]
%\item $A$ is a maximal uniform act,
%\item $A$ is an injective uniform act,
%\item $A$ is a minimal injective act,
%\item the endomorphism monoid of $A$ is a local monoid.
%\end{enumerate}
%\ete
Herein, we study conditions under which a uniform act is subdirectly irreducible.
\begin{definition}
 Over a semigroup $S$ an $S$-act $A$ is called $cocyclic$ if the intersection of its non-zero subacts is non-zero.
\end{definition}

Note that over a semigroup $S$ a right $S$-act is cocyclic if and only if it contains a least non-zero subact which is obviously simple or $\theta$-simple. Besides, every non-zero subdirectly irreducible $S$-act is cocyclic. For if $A$ is a subdirectly irreducible act and $\Sigma=\lbrace \rho_{_{B}}\mid \Theta\neq B\leqslant A\rbrace$ . Then $\bigcap \Sigma =\rho_C \neq\Delta_A$ where $C=\bigcap \limits _{ B\leqslant A \atop {\Theta\neq B}}B$ and  consequently $A$ is cocyclic.

Due to Proposition \ref{pr1}, Lemma \ref{le3} and the above argument the next theorem clarifies the structure of subdirectly irreducible acts over semigroups.
 \bte \lb{th1}
 For a semigroup $S$ a right $S$-act $A_S$ is subdirectly irreducible if and only if $A_S$ contains a large simple or $\theta$-simple subdirectly irreducible subact.
 \ete

%\bpf { \bf Necessity.} By assumption
%\[\bigcap \limits _ {{\substack {{B\leq A }\\ {{B\neq \Theta}}}}} {\rho}_{_{_{{B}}}} ={\rho}_{_{_{{(\bigcap \limits _ %{\substack {{B\leq A }\\ {{B\neq \Theta}}}}B)}}}}\neq \triangle.\]

%So $M_S=\bigcap \limits _ {{{\substack {{B\leq A }\\ {{B\neq \Theta}}}}}}B$ is a minimum nonzero large subact of $A_S$ %which is evidently simple or $\theta$-simple. Since any non-diagonal congruence on $M_S$ can be turned into a non-diagonal %one in $A_S$, $M_S$ is subdirectly irreducible.

%{\bf Sufficiency.} Let $M_S$ be a large simple or $\theta$-simple subdirectly irreducible subact of $A_S$. Since $M_S$ is large in %$A_S$, for any non-diagonal congruence $\rho \in {\rm Con}(A_S)$, $(\rho \, \cap \, M\times M)$ is a non-diagonal %congruence in $M_S$. Thus the following get the desired %result.
%\[\bigcap \limits _ {{\substack {{\rho \in {\rm Con}(A_S)}\\ {{\rho \neq \triangle}}}}} {\rho}\supseteq \bigcap \limits _ %{{\substack {{\rho \in {\rm Con}(A_S)}\\ {{\rho \neq \triangle}}}}} {(\rho \, {\cap \,M \times M)}}\neq \triangle.\]
%\epf\begin{corollary}
\bco \lb{co8}
Over a semigroup $S$, any uniform act with two zeros is subdirectly irreducible.
\eco
\bpf
Evidently the right $S$-act $\lbrace\theta_1, \theta_2\rbrace$ is a $\theta$-simple subdirectly irreducible act. Now Theorem \ref{th1} completes the proof.
\epf

\bco \lb{co1}
Let $A$ be an act over a semigroup $S$. Then $A$ is subdirectly irreducible if and only if  $A$ is  a uniform cocyclic $S$-act  which its ($\theta$-)simple subact is subdirectly irreducible.
\eco
\bpf
{\bf Necessity.} First note that  an act with a large simple or $\theta$-simple subact is cocyclic, indeed any large simple or $\theta$-simple subact is the least subact. Now by virtue of Theorem \ref{th1} and Proposition \ref{pr1} the result follows.

{ \bf Sufficiency.} Since $A$ is uniform its subdirectly irreducible ($\theta$-)simple subact is large and we are done in view of Theorem \ref{th1}.
\epf
Note that by Proposition \ref{pr1}, every irreducible act is uniform and so in the previous corollary we can replace ``uniform" with ``irreducible".
In light of Corollary \ref{co8}, every uniform act with two zeros is subdirectly irreducible. So in the next two results uniform acts with at most one zero element are characterized. The next lemma is an immediate result of the uniform notion.
\ble \lb{le1}
Let $S$ be a semigroup and $A$ be a non-zero  right $S$-act with at most one zero element. Then the following are equivalent:
\begin{enumerate}[{\rm i)}]
\item $A$ is uniform.
\item Every non-zero finitely generated subact of $A$ is a large subact.
\item Every non-zero cyclic subact of $A$ is a large subact.
\item Every non-zero indecomposable subact of $A$ is a large subact.
\item Every non-zero subact of $A$ is uniform.
\end{enumerate}
\ele
The equivalent condition (III) of the above lemma for uniform acts yields the following proposition which can be established by an adaptation of \cite[Theorem 7]{bert}.
\bpr \lb{pr8}
Let $A$ be a right act over a monoid $S$ with at most one zero element. $A$ is uniform if and only if for any non-zero element $a$ in $A$ and any two distinct elements $x$ and $y$ in $A$, there exist $s,t\in S$ such that $(as,at)\in \rho(x,y)$ (the monocyclic congruence generated by the pair $(x,y)$) and $(s,t)\notin{\rm ker}\lambda_a$.
\epr

\bco \lb{co3}
Let $S$ be an infinite cyclic group and $A$ be an $S$-act. Then $A$ is subdirectly irreducible if and only if $A$ is a finite uniform $S$-act which its simple subact is subdirectly irreducible.
\eco
\bpf
The necessity holds by Theorem 2.2.38 of \cite{kilp} and Corollary \ref{co1}. Also by Corollary \ref{co2}, we conclude that every uniform $S$-act with a finite number of subacts (in particular any finite uniform act) is cocyclic. Thus by Corollary \ref{co1} the result follows.
\epf
%\ble \lb{le6}
%Suppose $S$ is a semigroup and $A$ is a non-zero uniform $S$-act. If for a subact $B$ of $A$, $\theta_1\in B$ and $\theta_2\in A\backslash B$ are zero elements, then $B=\theta$.
%\ele
%\bpf
%Note that $\rho_B\cap \rho(\theta_1,\theta_2)=\Delta_A $. If $B\neq \theta$, then by Lemma 2.4, $\rho(\theta_1,\theta_2)=\Delta_A $ and so $\theta_1=\theta_2$, a contradiction.
%\epf
Recall that in a semigroup $S$ an element $s$ is called a left identity if $sx=x$ for any $x\in S$.

\bpr \lb{pr10}
Let $S$ be a semigroup and $A$ be a non-zero uniform right $S$-act. If $a\in A$ and for some left identity element $s \in S,as\neq a$, then the following hold:
\begin{enumerate} [{\rm i)}]
\item For every non-zero subact $B$ of $A$, $a\in B$.
\item $A$ is cocyclic.
\item $|Z(A)|\leq 1$.
\item If $b\in A$ and for some left identity element $t\in S, bt\neq b$, then $bS^{1}=aS^{1}$.
\end{enumerate}
\epr
\bpf
(i) Since $sx=x$ for every $x\in S$, it can be readily checked that $\rho(as,a)=\lbrace(as,a), (a,as)\rbrace \cup \Delta_A $. As $A$ is uniform any non-zero subact $B$ is large in $A$ and hence $\rho(as,a)\cap \rho_B\neq  \Delta_A $. Consequently $(a, as)\in \rho_B$ which gives $a\in B$.\\
(ii) This is a straightforward consequence of the first assertion.\\
(iii) Since $a\notin Z(A)$, in view of the first assertion $|Z(A)|\leq 1$.\\
(iv) Since $bS^1$ is a non-zero subact, the first assertion gives $a\in bS^{1}$. Since for some left identity element $t\in S,bt\neq b$ , applying the first assertion for $b$ implies that $b\in aS^1$.
\epf
The following is an immediate result of the above proposition.
%\bco \lb{co12}
%Let $A$ be a non-zero uniform $S$-act over a right zero semigroup $S$. Then $A$ is cocyclic.
%\eco
\bco \lb{co13}
Any non-zero uniform right act $A$ over a right zero semigroup $S$ has one of the following structures:
\begin{enumerate} [{\rm i)}]
\item $A=\theta\amalg \theta $,
\item $A=aS^{1}$,
\item $A=aS^{1}\amalg \theta$,\\
\end{enumerate}
which in the two latter cases $|Z(A)|\leq 1$ and $aS^1$ is simple or $\theta$-simple. Therefore any non-zero uniform act is cocyclic.
\eco
\bpf
Suppose that $A\neq \Theta\amalg\Theta$. Thus $A\backslash Z(A)\neq \emptyset$ and hence any non-zero element $a\in A$ satisfies the posed conditions in Proposition \ref{pr10}. Thus $|Z(A)|\leq 1$. On the other hand, for $a,b\in A\backslash Z(A)$, $as\neq a$ and $bt\neq b$ for some $s,t\in S$ and  Proposition \ref{pr10} yields $aS^{1}=bS^{1}$. Hence $A=aS^{1}$, if indecomposable, and $A=aS^{1}\amalg \Theta$ otherwise. Now the last statement is readily followed.
\epf
In \cite{kozh}, subdirectly irreducible acts over right zero semigroups are characterized. In what follows we do the same with another approach.

\bte \lb{th2}
Let $A$ be an act over a right zero semigroup $S$. $A$ is (subdirectly) irreducible if and only if $|A|=2$ or $A=B \sqcup \Theta$ where $B$ is a simple act of order $2$.
\ete
\bpf {\bf Necessity.} Let $A$ be an irreducible act with $|A|> 2$ over a right zero semigroup $S$. Using Corollary \ref{co13}, $|Z(A)|\leq 1$ and hence $|A\backslash Z(A)|\geq 2$. Now take two distinct elements $a$ and $b$ in $A \backslash Z(A)$. First note that if $as=\theta$ for some $s\in S$, then for any $t\in S$ $at=a(st)=(as)t=\theta t=\theta$. Now Proposition \ref{pr10}, part $(iv)$ implies that $bS^1=aS^1=\{a,\theta\}$, leading to a contradiction. Then $aS^1,bS^1\subseteq A\backslash Z(A)$. Let $as\neq a$ and $bt\neq b$ for some $s,t\in S$. Thus $\rho=\lbrace(as,a), (a,as)\rbrace \cup \Delta_A $ and $\sigma=\lbrace(bt,b), (b,bt)\rbrace \cup \Delta_A $ are two nondiagonal congruences over $A$. Our assumption provides that $\rho \cap \sigma \neq \Delta_A$. Since $a\neq b$, $(a,as)=(bt,b)$ and hence  $a=bt$ and $b=as$. Therefore $A \backslash Z(A)=aS^1=\{a,as\}$ and hence $A=aS^1 \sqcup \Theta$ where $aS^1$ is a simple act of order 2.

{\bf Sufficiency.} Suppose that $A=B \sqcup \Theta$ where $B$ is a simple act of order $2$. Then $\rho_B$ is the only proper congruence on $A$.
\epf

 \section{Uniform semigroups and injective properties on uniform acts}

 Herein we study interrelations among uniform notion and injective notions. First we mention necessary ingredients to start. Recall that an act $A$ is called C-injective (F-injective) if for any act $N$, any cyclic (finitely generated) subact $M$ of $N$, and any homomorphism $f\in$ Hom$(M,A)$, there exists a homomorphism $g\in$ Hom$(N,A)$ which extends $f,$ i.e., $g|_M=f$ (see \cite{Zhang}). We shall call an $S$-act $A$  \emph{indecomposable domain injective} (InD-injective for short) if it is injective relative to all inclusions from indecomposable acts.

Clearly InD-injective acts are C-injective. We refer the reader to \cite{MSN} for an expanded account on InD-injective acts over monoids.

\ble \lb{le7}
Let $A$ be an act over a monoid $S$. $A$ is C-injective if and only if $E(C)\subseteq A$, for any cyclic subact $C$ of $A$.
\ele

\bpf
Let $C$ be a cyclic subact of $A.$ Since $A$ is C-injective, the inclusion map from $C$ into $A$ can be extended to a homomorphism $g: E(C)\longrightarrow A.$ Moreover, since $C$ is large in $E(C)$, $g$ is a monomorphism. Thus $E(C) \subseteq A.$

For the converse, suppose that $D$ is a cyclic right $S-$act contained in a right $S-$act $B$ and $f: D\longrightarrow A$ is a homomorphism. Then $C=f(D)$ is a cyclic subact of $A$ and so by assumption $E(C)\subseteq A.$ Thus $f$ can be considered as a homomorphism from $D$ into $E(C).$ Hence there is an extension $\overline{f}$ of $f$ from $B$ into $E(C)$ which is also an extension into $A$. So $A$ is C-injective.
\epf

\bpr \lb{pr2}
Let $S$ be a monoid and $A$ be a non-trivial uniform right $S$-act. Then the following statements hold:
\begin{enumerate} [{\rm i)}]
\item $A$ is C-injective if and only if $A$ is injective.
\item If $B$ is a non-zero injective $S$-act, then every
monomorphism $f: B\longrightarrow A$ is an isomorphism.
\item If $B$ is a non-zero projective $S$-act, then every epimorphism $f: A\longrightarrow B$ is an isomorphism.
\end{enumerate}
\epr

\bpf
 Since $A$ is non-trivial uniform it contains a non-zero cyclic subact namely $C$. Now the first statement follows by Lemma \ref{le7} and Corollary \ref{co10}.\\
As an immediate consequence of the definition one may easily observe that any retract of a uniform act is isomorphic to itself. Now the stated conditions in parts ii and iii, provide that $B$ is a retract of the uniform act $A$ and hence $B$ is isomorphic to $A$.
\epf
We shall associate for each right $S$-act $A$ the right $S$-act $A^{\theta}_{S}$ with a zero element as below: \[ A^{\theta}_{S} = \begin{cases} A_S & \text{if } A_S \text{ contains a zero element,} \\ A_S \cup \Theta & \text{otherwise.}   \end{cases} \] In the next proposition we show that for uniform acts many of the known injective notions are equivalent.
\bpr \lb{pr3}
Let $S$ be a monoid and $A$ be a non-trivial uniform $S$-act. Then the following statements are equivalent:
\begin{enumerate} [{\rm i)}]
\item $A$ is C-injective.
\item $A$ is injective.
\item $A$ is InD-injective.
\item $A$ is F-injective.
\item $A^{\theta}_{S}$ is injective.
\end{enumerate}
\epr
\bpf
 Due to the implications ``Injective  $\Longrightarrow$ F-injective $\Longrightarrow$ C-injective", ``Injective  $\Longrightarrow$ InD-injective $\Longrightarrow$ C-injective" and the first part of Proposition \ref{pr2}, we reach to the equivalences of the parts i, ii, iii and iv. Moreover, parts iii and v are equivalent by \cite[Theorem 2.6]{MSN}.
\epf

Let $B$ be a large indecomposable subact of a decomposable act $A$. Then  $A=E \sqcup \Theta$ where $B\subseteq E$ and $E$ is indecomposable. For if $A=E \sqcup E'$ whence $E$ is an indecomposable component of $A$ containing $B$, then the canonical epimorphism $\pi :A\longrightarrow \ \dfrac{A}{E'}$ is a monomorphism since $\pi_{|B}$ is a monomorphism. Thus $E'=\Theta$.
%\ble \lb{le8}
%Let  $S$ be a semigroup and $A$ be a non-zero indecomposable $S$-act. If $E(A$ is decomposable, then $E(A= E_1\amalg \theta$,
%such that $A\subseteq E_1$.
%\ele
%\bpf
%Let $E(A)= E_1 \amalg E_2$ be a decomposition of $E(A)$. Since $A$ is indecomposable, either $A\subseteq E_1$ or $A\subseteq E_2$. Let $A \subseteq E_1$ and  Since $E_1\cap E_2=\emptyset$, . Since $A \subseteq {E_1} \subseteq E(A)$ and $A \subseteq' E(A)$, $ E_1 \subseteq'E(A$. Thus $\pi$ is a monomorphism and so $E_2=\theta$. Hence $E(A= E_1\amalg \theta$ and $A\subseteq E_1$.
%\epf

\bco \lb{co7}
Let $A$ be a non-zero indecomposable uniform $S$-act. Then $E(A)$ is indecomposable whenever $Z(A)\neq \emptyset$
\eco
\bpf
 If $E(A)$ is decomposable, the above argument implies that $E(A)= E\amalg \Theta$ such that $A \subseteq E$. Using Lemma \ref{le3}, $E(A)$ is uniform and since $E$ contains a zero element, Proposition \ref{pr9} implies that $E=\Theta$, a contradiction. Thus $E(A)$ is indecomposable.
 \epf
In what follows we investigate uniform acts in presence of a finiteness condition on congruences.
\begin{definition}
Let $S$ be a semigroup and $A$ be a right $S$-act. Then $A$ is called $strongly$ $right$ $noetherian$  if it satisfies the ascending  chain condition for right congruences.
\end{definition}
It is clear that every strongly right noetherian $S$-act is  right noetherian. It is known that for a noetherian module $M$ and a homomorphism $f: M \to M$, there exists a natural $n$ such that $\text{Im}f^n \cap \text{ker}f^n=0$, whence $f$ is an epimorphism if and only if it is an automorphism. In the next proposition we show the same is valid for acts over semigroups.

\bpr \lb{pr5} Let $A$ be a strongly right noetherian act over a semigroup $S$ and $f:A \to A$ be a homomorphism. Then ${\rm ker}(f^n) \cap \rho_{({\rm Im}f^n)}=\Delta_A$, for some natural $n$, whence $f$ is an epimorphism if and only if it is an isomorphism.
\epr
\bpf Suppose that $A$ is strongly right noetherian and $f:A \to A$ be a homomorphism. Regarding the ascending chain $\text{ker}(f) \subseteq \text{ker}(f^2) \subseteq \ldots$, there exists a natural $n$ such that the chain is stationary at $n$. Then $\text{ker}(f^n)=\text{ker}(f^{2n})$. Now if we take $B=\text{Im}(f^n)$ and $(a,b)\in \text{ker}(f^n) \cap \rho_B$ then $a=b$ or for some  $x,y \in A$, $a=f^{n}(x), b=f^{n}(y)$ and $f^{n}(a)=f^{n}(b)$. Thus $f^{n}(f^{n}(x))=f^{n}(f^{n}(y))$ and so $f^{2n}(x)=f^{2n}(y)$ which yields that $(x,y)\in \text{ker}(f^{2n})= \text{ker}(f^n)$. Therefore $a=f^{n}(x) =f^{n}(y)=b$ and then $\text{ker}(f^n) \cap \rho_B=\Delta_A$. Now if $f$ is an epimorphism then $\text{Im}f^n=A$ and regarding the condition ${\rm ker}(f^n) \cap \rho_{({\rm Im}f^n)}=\Delta_A$, $\text{ker}f^n$ is the diagonal relation and hence $f$ is a monomorphism as desired.
\epf
In module theory, if $M$ is a uniform noetherian module, then an endomorphism $f$ of $M$ is a monomorphism if and only if $f$ is not nilpotent. In the next two proposition we prove the same result for acts over semigroup, while in the first proposition the condition ``noetherian" is replaced with ``strongly noetherian".
\bpr \lb{pr11}
Suppose that $S$ is a semigroup and $A$ is a non-zero uniform and strongly right noetherian act. Then $f: A\longrightarrow A$ is a monomorphism if and only if it is not nilpotent (${\rm Im}f^n\neq \Theta$ for any natural $n$).
\epr
\bpf
We just need to prove the sufficiency. On account of Proposition \ref{pr5}, ${\rm ker}(f^n) \cap \rho_{({\rm Im}f^n)}=\Delta_A$, for some naturals $n$. Our assumption implies that $ \rho_{({\rm Im}f^n)}\neq \Delta_A$ and since $A$ is uniform ${\rm ker}(f) \subseteq {\rm ker}(f^n)= \Delta_A$ which completes the proof.
 \epf
\brm \lb{rm2}
It is worth to notice that if $A$ is  a uniform and strongly right noetherian right act with no zero element, then any endomorphism is a monomorphism. In particular, for any finite uniform act $A$ with no zero element,  every homomorphism $f:A\longrightarrow A$ is an isomorphism.
\erm
\bpr \lb{pr12}
Suppose that $S$ is a semigroup, $A$ is a non-zero noetherian uniform $S$-act and $f: A\longrightarrow A$ is a homomorphism. Then  $f^{-1}(\Theta)\neq \Theta $ if and only if $f^{m}=0$ for some natural $m$.
\epr
\bpf
{\bf Necessity.} Our assumption implies that for some $a\in A \backslash \Theta,f(a)= \theta \in \Theta$. For any natural number $n$, let $K_{n}=\lbrace{x\in A\mid f^{n}(x)=\theta}\rbrace$. By assumption $K_{1} $ is non-zero and hence $Z(A)\subseteq K_1$ considering the argument before Proposition \ref{pr9}. Since $f(Z(A))=\{\theta\}$, we have the ascending chain $K_{1}\subseteq K_{2}\subseteq \ldots$ which is stationary by assumption and then for some natural number $n$, $K_{n}=K_{n+1}=\cdots =K_{2n}$. If $x \in K_{n}\cap f^{n}(A)$, then $f^{n}(x)=\theta$ and $x=f^{n}(a)$ for some $a\in A$. Thus $\theta=f^{n}(x)=f^{n}(f^{n}(a))=f^{2n}(a)=f^{n}(a) =x$. Now since $ K_{n}\neq \Theta$, by Corollary \ref{co2}, $f^{n}(A)=\Theta$ and so $f^{n}=0$.

{\bf Sufficiency.} Suppose that $f^n=0$ and $f^{n-1}\neq 0$ for some natural $n$. Then for some $x\in A$, $f^{n-1}(x)\neq \theta$ and hence $f^{-1}(\Theta)\neq \Theta $
\epf

The rest of the paper is allocated to characterize some classes of uniform semigroups.
\bco \lb{co16}
Suppose $S$ is a uniform noetherian semigroup which contains at least a left zero element. If $x\in S$ and for a non-left zero element $y$, $xy$ is a left zero, then for some natural number $m$, $x^{m}$ is a left zero.
\eco
\bpf
In Proposition \ref{pr12}, let $f=\lambda_{x}$.
\epf

\bco \lb{co15}
If $S$ is a finite uniform semigroup which has no left zero element, then $S$ is left cancellative and right simple. Especially if $S$ is a monoid, then it is a group.
\eco
\bpf
In light of Remark \ref{rm2}, for every $x\in S$ the homomorphism $\lambda_{x}:S_S \longrightarrow S_S $ is an isomorphism and consequently $S$ is left cancellative and right simple. Now the second assertion follows immediately.
\epf
\ble
Let $S$ be a non-zero uniform semigroup and $e\in S$ be an idempotent element of $S$. Then $e$ is either a left zero or a left identity.
\ele
\bpf
If $e$ is not a left zero, then $eS$ is a large subact of $S_S$. Define $f:S_S\longrightarrow eS$ by $f(s)=es$ for every $s\in S$. Clearly $f$ is a well-defined homomorphism such that $f \mid_ {eS}$ is a monomorphism and hence $f$ is a monomorphism. Now, since for every $s\in S$ $f(s)=es=e(es)=f(es)$, we get  $s=es$, and therefore $e$ is a left identity.
\epf
\bpr
Let $S$ be a non-zero uniform monoid. Then $S=G\sqcup I$ where $G$ is the maximum subgroup of $S$ and $I$  is a two sided ideal of $S$.
\epr
\bpf
Let $G=\lbrace x\in S \mid  \exists\ y \in S, xy=1 \rbrace$. Then $G$ is a group. In fact if $xy=1$, then $yxyx=yx$. So $yx$ is an idempotent element of $S$ and by the previous lemma, $yx$ is a left zero or $yx=1$. If $yx$ is a left zero, then $x=xyx$ is a left zero and hence $xy=1$ implies that $x=1$, a contradiction. Thus $yx=1$. Also we can easily see that $ I=S\backslash G$ is a two sided ideal of $S$. Due to the definition of $G$, it contains any other subgroup and so is the maximum subgroup of $S$.
\epf
\bpr \lb{le4}
Let $S$ be a semigroup and $I$ be a non-zero large right ideal of $S$. If there exists an element $x \in S$, such that for any $m,n\in \mathbb{N},i,j\in I$, $x^mi\neq x^nj$ whenever $i\neq j$, in particular if $xi=i$ for any $i\in I$, then $x$ is a left identity in $S$.
\epr
\bpf
Let $\rho=\lbrace (a,b) \in S\times S\mid  \exists\ n, m \in  \mathbb{N}, x^{n}a=x^{m} b\rbrace$. Clearly $\rho $ is a right congruence on $S$. Our assumption implies that $\rho \cap \rho_I =\Delta _S$ and since $I$ is large in $S$, $\rho =\Delta_S$. Now suppose that $s \in S$ and $xs=z$. Then $x^{2} s= xz$ and hence $(s,z) \in \rho = \Delta_S$. Thus $xs=z=s$ and the result follows.
\epf
\bco \lb{co11}
Let $S$ be a uniform monoid and $x,y \in S$. Then  $xy=y$ implies that $x=1$ or $y$ is a left zero.
\eco

\bpf
If $xy=y$ and $y$ is not a left zero, then $yS$ is a large right ideal of $S_S$. Since $xy=y$, for every $t \in S$, $xyt=yt$ and by Proposition \ref{le4}, $xs=s$ for every $s \in S$. Hence for $s=1$ we obtain $x=1$.
\epf
\bco \lb{co17}
Suppose that $S$ is a uniform monoid and $G$ is the maximum subgroup. Then $G$ acts faithfully on non-zero elements, in particular, for every $s \in S \backslash Z(S_S)$, $| G|=| Gs|$.
\eco
\bpf
 Suppose $s\in S$ and define $f:G\longrightarrow Gs$ by $f(g)=gs$ for every $g \in G$. Now $f(g_1)=f(g_2)$ implies that $g_1s=g_2s$ and therefore $(g_2^{-1}g_1)s=s$. Thus Corollary \ref{co11} gives $g_2^{-1}g_1=1$ as desired.
\epf
In the next theorem we characterize regular monoids over which all cyclic acts are uniform.
\bte \lb{th3}
Let $S$ be a regular monoid. The following are equivalent:
\begin{itemize}
\item[{\rm (i)}]All cyclic $S$-acts are uniform.
\item[{\rm (ii)}] $S$ is uniform.
\item[{\rm (iii)}]$S$ meets one of the following three structures:
\begin{enumerate}
\item $S$ is a group,
\item $S=G\sqcup \{\theta\}$, whereas $G$ is a group and $\theta$ is a  zero.
\item $S=G\sqcup \{\theta_1,\theta_2\}$, whereas $G$ is a group, $\theta_1$ and $\theta_2$ are left zero elements, and $s\theta_i=\theta_j$ for all $s\in G \backslash \{1\}, 1\leq i\neq j \leq 2$.
\end{enumerate}
\end{itemize}
\ete
\bpf (i)$\longrightarrow$(ii). It is clear.

(ii)$\longrightarrow$(iii). Since $S$ is uniform, it has at most two left zeros. So in regard to this, three cases may occur.\\

{\bf $S$ has no left zero:}\\
Let $e$ be an idempotent in $S$. Since $e$ is not a left zero, and $ee=e$, Corollary \ref{co11} implies that $e=1$. Now, since $S$ is a regular monoid with only one idempotent, it is a group.\\

{\bf $S$ has only one left zero:}\\
Let $\theta$ be the left zero element of $S$ which is indeed the zero element and let $\theta \neq a \in S$. Let $x$ be an inverse of $a$. Since $axa=a$, $e=ax\neq \theta$, it follows that $e$ is an idempotent in $S$ which is not a left zero and hence $x$ is not a left zero. Analogously, $f=xa$ is an idempotent in $S$ which is not a left zero. Now Corollary \ref{co11} yields $e=1$, $f=1$ and $ax=xa=1$. Set $G=S \backslash \{\theta \}$. As any element in $G$ has an inverse, $G$ is closed under multiplication and hence is a group.

{\bf $S$ has two left zeros:}\\
Let $\theta_1$ and $\theta_2$ be two different left zero elements of $S$ and take $I=\{\theta_1,\theta_2\}$. Applying the same argument used in the second case, we conclude that $G=S\backslash \{\theta_1,\theta_2 \}$ is a group and $S$ can be presented of the form $S=G\sqcup I$. Now we turn to clarifying the multiplications of elements in $S$. To carry out this task we just need to consider the multiplications of the form $s\theta_i$ for all $s\in G \backslash \{1\}, 1\leq i \leq 2$. Indeed we prove that $s\theta_i=\theta_j$ for all $s\in G \backslash \{1\}, 1\leq i \leq 2$. Contrary to our claim, without restriction of generality, suppose that $s\theta_1=\theta_1$ for some $s\in G\backslash \{1\}$. Then $s^{-1}(s\theta_1)=(s^{-1}s)\theta_1=\theta_1=s^{-1}\theta_1$. This fact implies that $s\theta_2=\theta_2$. Indeed, if $s\theta_2=\theta_1$, then $s^{-1}(s\theta_2)=(s^{-1}s)\theta_2=\theta_2=s^{-1}\theta_1=\theta_1$, a contradiction. Now it can be readily checked that $(\theta_1, \theta_2)\notin \rho (s,1)$ therein. Now, $\rho(s,1) \cap \rho_I=\Delta_S$. Since $I$ is large, $\rho(s,1)=\Delta_S$ and hence $s=1$ which is leaded to a contradiction.

(iii)$\longrightarrow$(i). In the first case, since $S$ is a group, any cyclic act is simple and hence is uniform. In the second case, by virtue of the monoid structure, any cyclic act has no proper subact. In the last case, suppose that $aS$ is a cyclic act over $S$. Clearly, by virtue of the monoid structure, the only proper subact of $aS$ is $\{a\theta_1,a\theta_2\}$, which we prove that it is large in $aS$. For this purpose, suppose that $f$ is a homomorphism from $aS$ to an arbitrary act which is not monomorphism. To complete the proof we need only to verify that $f(a\theta_1)=f(a\theta_2)$. Let $f(as)=f(at)$ for two different elements $s,t \in S$. If $\{s,t\}=\{\theta_1,\theta_2\}$ then we are done. Otherwise, $\{s,t\}\cap G \neq \emptyset$ which leads to $f(a)=f(au)$ for some $u\neq 1$ and hence, $f(a\theta_1)=f(a\theta_2)$.\\
%In the case that $S$ is a periodic monoid, regarding the fact that each element has an idempotent power, using  an analogous argument as above, we can complete the proof
\epf

\brm \lb{re1}
Note that in the previous theorem, $S$ should be a monoid. Otherwise, for example if $S$ is a right zero semigroup, then $S$ is uniform which has none of the above structures.
\erm
 %The next proposition asserts a critical property of uniform acts which leads us to considerable results.

% Let $S$ be a monoid and $B_S$ be a non-zero subact of a uniform act $A_S$. Then for  and $a,b\in A$. Then {\rm Ker}$\lambda_a$ and {\rm Ker}$\lambda_b$ are not comparable, consequently, the mapping $F:A \to {\rm Con}(S)$ given by $F(a)={\rm Ker}\lambda_a$ is injective.
 In the next proposition we characterize monoids over which there exists a non-trivial uniform act with two zeros.
\bpr \lb{pr7} Let $S$ be a monoid. There exists a non-trivial uniform (subdirectly irreducible) act with two zeros if and only if $S$ is not left reversible.
\epr
\bpf {\bf Necessity.} Let $A_S$ be a non-trivial uniform act with two zeros $\theta_1$ and $\theta_2$. Let $a \in A \backslash \{\theta_1, \theta_2\}$. Since $aS$ is a non-zero subact of $A$, Corollary \ref{co2} implies that $\theta_1, \theta_2 \in aS$. Now the right ideals $I=\{s\in S\,|\, as=\theta_1\}$ and $J=\{s\in S\,|\, as=\theta_2\}$ have empty intersection and we are done.

{ \bf Sufficiency.} Suppose that $aS\cap bS=\emptyset$ for some $a,b\in S$. Since $1\notin aS\cup bS$, the Rees factor act $S/(aS\cup bS)$ is a right $S$ act with two different zeros namely $[a]$ and $[b]$. % and $|S/(aS\cup bS)\geq 3$. Let $\Sigma =\{(A,B)\,|\,A\text{~and~}B \text{~are right ideals of~}S,a\in A,\, b\in B \text{~and~} A\cap B=\emptyset\}$ be ordered by componentwise inclusion order. Let $\{(A_i,B_i)\,|\,i\in I\}$ be a chain in $\Sigma$. Evidently, $(\bigcup\limits_{i\in I}A_i,\bigcup\limits_{i\in I}B_i)$ is an upper bound of this chain in $\Sigma$.
Let $\Sigma$ be the set of right congruences $\rho$ on $S$ in which $[a]_{\rho}$ and $[b]_{\rho}$ are two different zero elements of $S/\rho$ that is ordered by inclusion. Then for any chain $(\rho_i)_I$ in $\Sigma$, $\bigcup \limits_I \rho_i$ is an upper bound in $\Sigma$. Therefore, Zorn's lemma yields a maximal congruence $\rho$ in $\Sigma$. Now $S/\rho$ is a non-trivial subdirectly  irreducible uniform act with two different
zeros. Indeed, any nondiagonal right congruence on $S/\rho$ connects the zero elements $[a]_{\rho}$ and $[b]_{\rho}$.\epf

\brm \lb{re2} Let $A_S$ be a uniform act with two zeros namely $\theta_1, \theta_2$. In light of Corollary \ref{co10}, $E(A)=E(\{\theta_1, \theta_2\})$. On the other hand, $\{0,1\}^S$ is an injective act with two zeros and hence containing $E(A)$. Therefore any uniform act with two zeros can be embedded in $\{0,1\}^S$. Thereby, we get an upper bound on the cardinality of uniform acts with two zeros.
\erm
\bco \lb{co9} Let $S$ be a semigroup. If $A_S$ is a uniform act with two zeros then $|A|\leq 2^{|S|}$.
\eco

%\textbf{Acknowledgment}
%The authors would like to thank the
%referee for providing valuable comments and suggestions.

\end{document}